\documentclass[12pt]{amsart}
\usepackage[latin1]{inputenc}
\usepackage{amscd}
\usepackage{latexsym}
\usepackage{pstcol,pst-plot,pst-3d}
\usepackage{mathptmx}
\usepackage{multicol}

\psset{unit=0.7cm,linewidth=0.8pt,arrowsize=2.5pt 4}

\newpsstyle{fatline}{linewidth=1.5pt}
\newpsstyle{fyp}{fillstyle=solid,fillcolor=verylight}
\definecolor{verylight}{gray}{0.97}
\definecolor{light}{gray}{0.9}
\definecolor{medium}{gray}{0.85}

\def\NZQ{\Bbb}               

\def\CC{{\NZQ C}}

\def\FFF{{\NZQ F}}
\def\TT{{\NZQ T}}
\def\opn#1#2{\def#1{\operatorname{#2}}} 
\opn\chara{char} \opn\length{\ell} \opn\pd{pd} \opn\rk{rk}
\opn\projdim{proj\,dim} \opn\injdim{inj\,dim} \opn\rank{rank}
\opn\depth{depth} \opn\grade{grade} \opn\height{height}
\opn\embdim{emb\,dim} \opn\codim{codim}

\opn\Tr{Tr} \opn\bigrank{big\,rank}
\opn\superheight{superheight}\opn\lcm{lcm}
\opn\trdeg{tr\,deg}%
\opn\reg{reg} \opn\lreg{lreg} \opn\skel{skel}
\opn\multideg{multideg}
\opn\div{div} \opn\Div{Div} \opn\cl{cl} \opn\Cl{Cl}
\opn\Spec{Spec} \opn\Supp{Supp} \opn\supp{supp} \opn\Sing{Sing}
\opn\Ass{Ass}
\opn\Ann{Ann} \opn\Rad{Rad} \opn\Soc{Soc}
\opn\Ker{Ker} \opn\Coker{Coker} \opn\Im{Im} \opn\Hom{Hom}
\opn\Tor{Tor} \opn\Ext{Ext} \opn\End{End} \opn\Aut{Aut}
\opn\id{id}

\opn\nat{nat}
\opn\pff{pf}
\opn\Pf{Pf} \opn\GL{GL} \opn\SL{SL} \opn\mod{mod} \opn\ord{ord}
\opn\aff{aff} \opn\con{conv} \opn\relint{relint} \opn\st{st}
\opn\lk{lk} \opn\cn{cn} \opn\core{core} \opn\vol{vol}
\opn\link{link} \opn\star{star} \opn\skel{skel}
\opn\gr{gr}

\def\pot#1#2{#1[\kern-0.28ex[#2]\kern-0.28ex]}

\opn\dirlim{\underrightarrow{\lim}}
\opn\inivlim{\underleftarrow{\lim}}
\let\union=\cup
\let\sect=\cap
\let\dirsum=\oplus

\let\Union=\bigcup
\let\Sect=\bigcap
\let\Dirsum=\bigoplus

\let\to=\rightarrow
\let\To=\longrightarrow
\def\Implies{\ifmmode\Longrightarrow \else
     \unskip${}\Longrightarrow{}$\ignorespaces\fi}
\def\implies{\ifmmode\Rightarrow \else
     \unskip${}\Rightarrow{}$\ignorespaces\fi}
\def\iff{\ifmmode\Longleftrightarrow \else
     \unskip${}\Longleftrightarrow{}$\ignorespaces\fi}

\let\:=\colon
\newtheorem{Theorem}{Theorem}[section]
\newtheorem{Lemma}[Theorem]{Lemma}
\newtheorem{Corollary}[Theorem]{Corollary}
\newtheorem{Proposition}[Theorem]{Proposition}

\let\epsilon\varepsilon
\let\phi=\varphi
\let\kappa=\varkappa
\def\qed{\ifhmode\textqed\fi
   \ifmmode\ifinner\quad\qedsymbol\else\dispqed\fi\fi}
\def\textqed{\unskip\nobreak\penalty50
    \hskip2em\hbox{}\nobreak\hfil\qedsymbol
    \parfillskip=0pt \finalhyphendemerits=0}
\def\dispqed{\rlap{\qquad\qedsymbol}}

\def\KK{{\mathcal K}}
\def\BB{{\mathcal B}}

\def\LL{{\mathcal L}}
\def\II{{\mathcal I}}
\def\JJ{{\mathcal J}}

\opn\inii{in} \opn\inim{inm} \opn\rate{rate}
\begin{document}
\title{The monomial ideal of a finite meet-semilattice}
\author{J\"urgen Herzog, Takayuki Hibi and Xinxian Zheng}
\address{J\"urgen Herzog, Fachbereich Mathematik und
Informatik, Universit\"at Duisburg-Essen, 45117 Essen, Germany}
\email{juergen.herzog@uni-essen.de}
\address{Takayuki Hibi, Department of Pure and Applied Mathematics,
Graduate School of Information Science and Technology, Osaka
University, Toyonaka, Osaka 560-0043, Japan}
\email{hibi@math.sci.osaka-u.ac.jp}
\address{Xinxian Zheng, Fachbereich Mathematik und
Informatik, Universit\"at Duisburg-Essen, 45117 Essen, Germany}
\email{xinxian.zheng@uni-essen.de}
\date{}
\subjclass{13D02, 13H10, 06A12, 06D99} \maketitle

\begin{abstract}
Squarefree monomial ideals  arising from finite meet-semilattices
and their free resolutions are studied. For the squarefree
monomial ideals corresponding to poset ideals in a distributive
lattice the Alexander dual is computed.
\end{abstract}

\section*{Introduction}
One of the most influential results in the classical lattice
theory is Birkhoff's fundamental structure theorem for finite
distributive lattices \cite[Theorem 3.4.1]{St}, which guarantees
that, given a finite distributive lattice $\LL$, there is a unique
poset (partially ordered set) $P$ such that $\LL$ is isomorphic to
the poset $\JJ(P)$ consisting of all poset ideals (including the
empty set) of $P$, ordered by inclusion.  (A poset ideal of $P$ is
a subset $I \subset P$ with the property that if $p \in I$ and $q
\in P$ with $q \leq p$, then $q \in I$.) In fact, if $P$ is the
subposet of $\LL$ consisting of all join-irreducible elements of
$\LL$, then $\LL = \JJ(P)$. (An element $p \in \LL$ with $p \neq
\hat{0}$ is called join-irreducible if there is no $q, r \in \LL$
with $q < p$ and $r < p$ such that $p = q \vee r$.) In other
words, by identifying $\LL$ with $\JJ(P)$, if $p \in \LL$ and $I =
\{ q \in P \: q \leq p \} \in \JJ(P)$, then $p = I$.

Fix a finite distributive lattice $\LL = \JJ(P)$. Let $K$ be a
field and $S = K[\{ x_p, y_p \}_{p \in P}]$ the polynomial ring in
$2|P|$ variables over $K$ with  $\deg x_p = 1$ and $\deg y_p=1$
for all $p\in P$. We associate each element $I \subset \JJ(P) =
\LL$ with the squarefree monomial $u_I = (\prod_{p \in I}x_p)
(\prod_{p \in P \setminus I}y_p) \in S$.  In the previous paper
\cite{HH} the monomial ideal $H_\LL = (u_I)_{I \in \LL}$ is
discussed from viewpoints of both combinatorics and commutative
algebra. The purpose of the present paper is to introduce the
squarefree monomial ideal $H_\LL$ for an arbitrary finite
meet-semilattice $\LL$ and to generalize some of the results
obtained in \cite{HH}.

Now, let $\LL$ be an arbitrary finite  meet-semilattice \cite[p.
103]{St} and $P \subset \LL$ the set of join-irreducible elements
of $\LL$. For each element $q \in \LL$ we write $\ell(q) = \{ p
\in P \: p \leq  q \} \subset P$.  In particular $\ell(\hat{0}) =
\emptyset$.  Note that $\ell(q)$ is a poset ideal of $P$, and that
$q \in \ell(q)$ if and only if $q$ is join-irreducible.  We thus
obtain the map $\ell \: \LL \to \BB_P$, which we call the
canonical embedding of $\LL$ into the Boolean lattice $\BB_P$
consisting of all subsets of $P$ ordered by inclusion.  As in the
case of finite distributive lattices explained in the previous
paragraph, let $K$ be a field and $S = K[\{ x_p, y_p \}_{p \in
P}]$ the polynomial ring in $2|P|$ variables over $K$ with  $\deg
x_p = 1$ and $\deg y_p=1$ for all $p\in P$.  We associate each
element $q \in \LL$ with the squarefree monomial $u_q = (\prod_{p
\in \ell(q)} x_p) (\prod_{p \in P \setminus \ell(q) }y_p) \in S$
and set $H_\LL = (u_q)_{q \in P} \subset S$.

In the present paper the following topics on squarefree monomial
ideals $H_\LL$ arising from finite meet-semilattices $\LL$ will be
studied:
\begin{itemize}
\item
When has the squarefree monomial ideal $H_\LL$
a linear resolution?  Theorem 1.3 guarantees
that $H_\LL$ has a linear resolution if and only if
$\LL$ is meet-distributive.
(A finite meet-semilattice $\LL$ is called
meet-distributive if each interval
$[x,y]=\{p\in\LL\: x\leq p\leq y\}$ of $\LL$ such that
$x$ is the meet of the lower
neighbors of $y$ in this interval is Boolean.
Here we call $z$ a lower neighbor of $y$ if $y$ covers $z$.)
\medskip
\item How can we construct a finite multigraded free
$S$-resolution $\FFF$ of $H_\LL$?  A construction of such a finite
free resolution is given in Theorem \ref{free} (a). Moreover, we
will characterize when our resolution is minimal.  In fact, it
will be proved in Theorem \ref{free} (b) that our resolution is
minimal if and only if, for any $p \in \LL$ and for any proper
subset $S \subset N(p)$ the meet $\bigwedge\{q\: q\in S\}$ is
strictly greater than the meet $\bigwedge\{q\:q\in N(p)\}$, where
$N(p)$ is the set of lower neighbors of $p$ in $\LL$. In
particular, if $\LL$ is a meet-distributive meet-semilattice, then
our finite free resolution is minimal (Corollary \ref{meet}). On
the other hand, when $\LL$ is a meet-distributive
meet-semilattice, the differential $\partial$ in the finite
multigraded free $S$-resolution $\FFF$ of $H_\LL$ obtained in
Theorem \ref{free} (a) will be described (Theorem \ref{maps}).
\medskip
\item Since $H_\LL$ is a squarefree monomial ideal, there is a
simplicial complex $\Delta$ whose Stanley--Reisner ideal
$I_\Delta$ coincides with $H_\LL$. We are interested in the
Alexander dual $\Delta^\vee$ of $\Delta$.  In case that $\LL$ is a
finite distributive lattice, a nice description of $\Delta^\vee$
can be obtained (\cite[Lemma 3.1]{HH}). It seems, however, rather
difficult, for an arbitrary finite meet-semilattice, to obtain an
explicit description of the Alexander dual of $H_\LL$. We will
consider a special meet-distributive meet-semilattice, namely, a
poset ideal $\II$ of a finite distributive lattice.  In this case
a combinatorial description of the Alexander dual of $H_\II$ can
be obtained (Theorem 4.2). Moreover, since $H_\II$ has a linear
resolution, it follows that the Alexander dual of $H_\II$ is
Cohen--Macaulay.  The combinatorics on such  Cohen--Macaulay
complexes is discussed in Theorem \ref{takayuki}.
\end{itemize}

\section{Algebraic characterizations of meet-distributive
meet-semilattices}

Let $\LL$ be an arbitrary finite meet-semilattice (c.f.\ \cite[p.
103]{St}), and $P\subset \LL$ the set of join-irreducible elements
of $\LL$. We denote by $\hat{0}$ and $\hat{1}$ the minimal and
maximal element of $\LL$. (Since $\LL$ is a finite
meet-semilattice, it follows \cite[Proposition 3.3.1]{St} that
$\LL$ possesses $\hat{1}$ if and only if $\LL$ is a lattice.)
Recall that $p\in \LL$ is {\em join-irreducible} if $p\neq
\hat{0}$ and $p$ is not a join of elements strictly less than $p$.

To each element $p\in \LL$ we associate the subset $\ell(p)=\{q\in
P\: q\leq p\}$ of $P$. Note that $p\in \ell(p)$ if  and only if
$p$ is join irreducible. In any case, $\ell(p)$ is a poset ideal
of $P$. Recall that a {\em poset ideal} of $P$ is a subset
$I\subset P$ such that if $r\in I$ and $t\leq r$, then $t\in I$.
The {\em set of generators} of $I$ is the set of maximal elements
in $I$, denoted by $G(I)$.

We obtain  a map
\[
\ell\: \LL\To \BB_P,
\]
which we call the {\em canonical embedding} into the Boolean
lattice $\BB_P$ consisting of all subsets of $P$ ordered by
inclusion.

We call the cardinality of $\ell(p)$ the {\em degree} of $p$, and
denote it by $\deg p$.  One always has the inequality $\rank p\leq
\deg p$. Recall that the {\em rank} of $p$ is the maximal length
of chains descending from $p$.

\begin{Lemma}
\label{better} Let ${\mathcal L}$ be a finite meet-semilattice,
$\ell$ the canonical embedding and  $s,t\in {\mathcal L}$ any two
elements. We have
\begin{enumerate}
\item[(i)]  $s=t$ if and only if $\ell(s)= \ell(t)$;

\item[(ii)] $s\leq t$ if and only if $\ell(s)\subseteq \ell(t)$;

\item[(iii)] $\ell(s)\cap \ell(t)=\ell(s\wedge t)$.
\end{enumerate}
\end{Lemma}
\begin{proof}
Note that each element of $\LL$ is the join of elements in $P$.
 From this observation all assertions follow.
\end{proof}

The lemma implies that $\ell$ is an injective order preserving
map. In general however, $\ell$ is not an embedding of lattices.
It is not difficult to see that $\ell$ is an embedding of
meet-semilattices if and only if $\LL$  is meet-distributive.

\medskip
We now introduce the definition of meet-distributive
meet-semilattices. A finite meet-semilattice $\LL$ is called {\em
meet-distributive} if  each interval $[x,y]=\{p\in\LL\: x\leq
p\leq y\}$ of $\LL$ such that $x$  is the meet of the lower
neighbors of $y$ in this interval is Boolean. Here we call $z$ a
lower neighbor of $y$ if $y$ covers $z$.

The following combinatorial characterization of meet-distributive
lattices are discussed in the survey article \cite{Ed}. A finite
meet-semilattice is called {\em graded} if for each  elements all
of its maximal chains have the same length.

\begin{Lemma}
\label{edelman} For a finite lattice $\LL$ the following
conditions are equivalent:
\begin{enumerate}
\item[(i)] $\LL$ is meet-distributive; \item[(ii)] $\LL$ is graded
and $\deg \hat{1} =\rank \hat{1}$; \item[(iii)] $\LL$ is graded
and $\deg \hat{p} =\rank \hat{p}$ for all $p\in\LL$; \item[(iv)]
each element in $\LL$ is a unique minimal join of join-irreducible
elements; \item[(v)] $\LL$ is lower semimodular, and any upper
semimodular sublattice is distributive.
\end{enumerate}
\end{Lemma}

We now introduce   the squarefree monomial ideal $H_\LL$
associated with a finite meet-semilattice $\LL$. Let $P$ be the
set of join irreducible elements of $\LL$. Let $K$ be a field and
$S=K[\{x_p,y_p\}_{p\in P}]$ the polynomial ring in $2|P|$
variables over $K$. For each element $q\in\LL$ write
\[
u_q=\prod_{p\in \ell(q)}x_p\prod_{p\in P\setminus \ell(q)}y_p,
\]
and set $H_\LL=(u_q)_{q\in\LL}$.

Note that $\height(H_\LL)=2$ if $\LL$ is a lattice. In fact,
$H_\LL\subset (x_p,y_p)$ for any $p\in P$ while on the other hand
$u_{\hat{0}}= \prod_{p\in P}y_p$ and $u_{\hat{1}}=\prod_{p\in
P}x_p$ both belong to $H_\LL$ and have no common factor.

 Let $I$ be a monomial ideal with the (unique) minimal set $G(I)$ of
monomial generators. The ideal $I$ is said to have {\em linear
quotients} if the elements of $G(I)$ can be ordered $u_1, \ldots,
u_m$ such that the colon ideals $(u_1,\ldots, u_{i-1}):u_{i}$ are
generated by variables. If $I$ is squarefree, then $I$ has linear
quotients if and only if for each $i$ and each  $j<i$  there
exists $k<i$ such that $u_k/[u_k,u_i]$ is a variable and divides
$u_j$. Here $[u,v]$ denotes the greatest common divisor of $u$ and
$v$.

It is easy to see that if all generators of $I$ have the same
degree, and  $I$ has linear quotients, then $I$ has a linear
resolution.

We now come to our algebraic characterization of meet-distributive
meet-semilattices.

\begin{Theorem}
\label{main} Let $\LL$ be an arbitrary finite meet-semilattice.
The following conditions are equivalent:
\begin{enumerate}
\item[(i)] $\LL$ is meet-distributive; \item[(ii)] $H_\LL$ has
linear quotients; \item[(iii)] $H_\LL$ has a linear resolution;
\item[(iv)]  $H_\LL$ has linear relations.
\end{enumerate}
\end{Theorem}

\begin{proof}
(i)\implies (ii): We fix a linear order $\prec$ on $\LL$ which
extends the partial order given by the degree. We put $u_r<u_q$ if
$r\prec q$. For any $u_q\in H_\LL$ and any $u_r < u_q$, let $t$ be
a lower neighbor of $q$ in the interval $[r\wedge q, q]$. Then
$u_t/[u_t,u_q]=y_p$, where $\{p\}=\ell(q)\setminus\ell(t)$. We
claim that $y_p$ divides $u_r$. If not, then $x_p$ divides $u_r$
and so $p\in \ell(r)\sect\ell(q)=\ell(r\wedge q)$. Thus
$p\in\ell(t)$, since $r\wedge q\leq t$.

(ii)\implies (iii) and (iii)\implies (iv) are trivial.

(iv) \implies (i): Suppose $\LL$ is not meet-distributive. Then by
Lemma \ref{edelman}(iii) (which is also valid if $\LL$ is a
meet-distributive meet-semilattice) there exist $p,q\in\LL$ such
that $q$ is lower neighbor of $p$ and $\deg p-\deg q>1$. The ideal
$(u_p,u_q)$ is generated by precisely those monomials in
$G(H_\LL)$ which are not divided by $x_r$ for all $r\in P\setminus
\ell(p)$, and are not divided by all $y_s$ for all $s\in\ell(q)$.
Since we assume that $H_\LL$ has linear relations, the restriction
lemma in \cite[Lemmma 4.4]{HHZ} implies that $(u_p,u_q)$ has
linear relations contradicting the fact that $\deg p-\deg q>1$.
\end{proof}

\begin{Corollary}
\label{upper semi} Let $\LL$ be a finite upper semimodular
lattice. Then the following conditions are equivalent:
\begin{enumerate}
\item[(i)] $H_\LL$ has a linear resolution; \item[(ii)] $\LL$ is
distributive.
\end{enumerate}
\end{Corollary}

\begin{proof}
The assertion follows from Lemma \ref{edelman}(v) and Theorem
\ref{main}.
\end{proof}

Let $\Delta$ be a simplicial complex on the vertex set
$[n]=\{1,\ldots,n\}$. The simplicial complex
\[
\Delta^\vee=\{[n]\setminus F\: F\not\in\Delta\}
\]
is called the {\em Alexander dual}  of $\Delta$. It is easy to see
that $(\Delta^\vee)^\vee=\Delta$.

A {\em vertex cover} of $\Delta$ is a set $G\subset [n]$ such that
$G\sect F \neq\emptyset$ for all $F\in {\mathcal F}(\Delta)$ where
${\mathcal F}(\Delta)$ denotes the set of facets (maximal faces)
of $\Delta$. A vertex cover is called {\em minimal} if it is
minimal with respect to inclusion. We also \texttt{}denote by
${\mathcal C}(\Delta)$ the set of minimal vertex covers of
$\Delta$.

As usual we denote by $I_\Delta$ the Stanley--Reisner ideal of
$\Delta$. The {\em facet ideal} is defined to be
\[
I(\Delta)=(x_F\: F\in{\mathcal F}(\Delta)),
\]
where $x_F=\prod_{i\in F}x_i$.

For $F=\{i_1,\ldots, i_k\}\subset [n]$ set  $P_F=(x_{i_1},\ldots,
x_{i_k})$, and let $\Gamma$ be the unique  simplicial complex such
that $I_\Delta=I(\Gamma)$. Then
\begin{eqnarray}
\label{cover} I_{\Delta}=\Sect_{F\in {\mathcal C}(\Gamma)}
P_F\quad \text{and} \quad I_{\Delta^\vee} = (x_F\: F\in {\mathcal
C}(\Gamma)).
\end{eqnarray}

Set $F^c=[n]\setminus F$, and $$\Delta^c=\langle F^c\: F\in
{\mathcal F}(\Delta)\rangle.$$ Then \
\begin{eqnarray}
\label{formula1}
 I_{\Delta^\vee}=I(\Delta^c).
\end{eqnarray}

The following lemma gives important algebraic properties of
Alexander duality.

\begin{Lemma}
\label{duality} Let $K$ be a field, $\Delta$ a simplicial complex,
$I_\Delta$ the Stanley--Reisner ideal and $K[\Delta]$ the
Stanley--Reisner ring of $\Delta$. Then
\begin{enumerate}
\item[(i)] {\em (Eagon--Reiner \cite{ER})} $K[\Delta]$ is
Cohen-Macaulay \iff
 $I_{\Delta^\vee}$ has a linear resolution.
\item[(ii)] {\em (\cite{HHZ})} $\Delta$ is shellable \iff
$I_{\Delta^\vee}$ has linear quotients.
\end{enumerate}
\end{Lemma}

Theorem \ref{main} together with Lemma \ref{duality} yields

\begin{Corollary}
\label{clear} Let $\LL$ be an arbitrary finite meet-semilattice,
and let $\Delta_\LL$ be the simplicial complex whose
Stanley--Reisner ideal is $H_\LL$. The following conditions are
equivalent:
\begin{enumerate}
\item[(i)] $(\Delta_\LL)^\vee$ is shellable; \item[(ii)]
$(\Delta_\LL)^\vee$ is Cohen--Macaulay; \item[(iii)] $\LL$ is
meet-distributive.
\end{enumerate}
\end{Corollary}

\begin{Proposition}
\label{minimal} Let $\LL$ be  a finite lattice and $P$ its poset
of join irreducible elements. Then
\begin{enumerate}
\item[(i)] the minimal prime ideals of height $2$ of $H_\LL$ are
$(x_p,y_q)$ where $p,q\in P$ and $p\leq q$; \item[(ii)] $H_\LL$
has only height 2 minimal prime ideals if and only if
 $\LL$ is distributive.
\end{enumerate}
\end{Proposition}

\begin{proof}
Let $\hat{\LL}$ be the distributive lattice consisting of all
poset ideals of $P$. Then $\ell$ induces an injective order
preserving map $\ell\:\LL\to \hat{\LL}$. Thus $H_\LL\subset
H_{\hat{\LL}}$, and equality holds if and only if $\LL$ is
distributive. This follows from Birkhoff's   fundamental structure
theorem \cite{St}.

(i) The minimal prime ideals of $H_{\hat{\LL}}$ are precisely the
ideals $(x_p,y_q)$ where $p,q\in P$ and $p\leq q$, see \cite{HH}.
Of course these are also minimal prime ideals of $H_\LL$. We claim
that there are no other minimal prime ideals of height 2 of
$H_\LL$. Indeed, any such ideal must contain some $x_p$ and some
$y_q$, since $\prod_{p\in P}x_p$ and $\prod_{p\in P}y_p$ belong to
$H_\LL$.  Suppose $p\not\leq q$, then $u_q$ is not contained in
$(x_p, y_q)$.

(ii) It remains to show that if $\LL$ is not distributive, then
there exists a minimal prime ideal of $H_\LL$ of height $>2$. In
fact, the proof of (i) shows that if such a minimal prime ideal
does not exist, then $H_\LL=H_{\hat{\LL}}$. Therefore
$\LL=\hat{\LL}$, and hence $\LL$ is distributive.
\end{proof}

Proposition \ref{minimal} together with (\ref{cover}) implies

\begin{Corollary}
\label{flag} A finite lattice $\LL$ is distributive if and only if
$(\Delta_\LL)^\vee$ is flag.
\end{Corollary}

\section{A free resolution of $H_\LL$}

The main theorem of the present section is the following

 \begin{Theorem}
 \label{free} Let $\LL$ be finite meet-semilattice.
 \begin{enumerate}
 \item[(a)]
  There exists a finite multigraded free $S$-resolution $\FFF$ of
$H_\LL$ such that for each $i\geq 0$, the free module $F_i$ has a
basis with basis elements
\[b(p;S)\]
 where $p\in\LL$ and $S$  is a subset of the set of lower neighbors
$N(p)$ of  $p$ with  $|S|=i$.

 The multidegree of $b(p;S)$ is the least common multiple of $u_p$ and
all monomials $u_q$ with $q\in S$.

\item [(b)] The following conditions are equivalent:
\begin{enumerate}
\item [(i)] the resolution constructed in {\em (a)}  is minimal;
\item[(ii)] for any $p\in\LL$ and  any proper subset $S\subset
N(p)$ the meet $\bigwedge\{q\: q\in S\}$ is strictly greater than
the meet $\bigwedge\{q\:q\in N(p)\}$.
\end{enumerate}
\end{enumerate}
 \end{Theorem}

We call a finite meet-semilattice satisfying condition (b)(ii)
{\em meet-irredundant}.

 \begin{proof}[Proof of \ref{free}]
(a) The resolution will be built by an iterated mapping cone
construction. As in the proof of Theorem \ref{main} we fix a
linear order $\prec$ on $\LL$ which extends the partial order
given by the degree. For any $p$ in $\LL$ we construct inductively
a complex $\FFF(p)$ which is a multigraded free $S$-resolution of
the ideal $H_\LL(p)$ generated by all $u_q\in H_\LL$ with
$q\preceq p$. Then $\FFF(q)$ is the desired resolution, where
$q\in\LL$ is the maximal element with respect to $\prec$.

The complex $\FFF(\hat{0})$ is defined as $F_i(\hat{0})=0$ for
$i>0$, and $F_0(\hat{0})=S$. This complex together with the
augmentation map $\epsilon\: S\to H_\LL(\hat{0})$, $1\mapsto
u_{\hat{0}}$ is a free resolution of $H_\LL(\hat{0})$.

Now let $p\in \LL$, $p\neq \hat{0}$, and let $q\in \LL$, $q\prec
p$ be the element preceding $p$. Then $H_\LL(p)=(H_\LL(q),u_p)$,
and hence we get an exact sequence of multigraded $S$-modules
\[
0\To (S/L)(-\multideg u_p)\To S/H_\LL(q)\To S/H_\LL(p)\To 0,
\]
where $L$ is the colon ideal $H_\LL(q):u_p$. As in the proof of
\ref{main} one shows that
\[
L=(\{u_t/[u_t,u_p]\}_{t\in N(p)}).
\]
Let $\TT$ be the Taylor complex associated with the monomials
$u_t/[u_t,u_p]$, $t\in N(p)$, see \cite{Ei}. Then $\TT$ is a
multigraded free resolution of $S/L$ with $T_0=S$,
$T_1=\Dirsum_{t\in N(p)} Se_t$ and $T_i=\bigwedge^iT_0$ for $i\geq
1$. Thus $T_i$ has a basis whose elements are $e_{t_1}\wedge
e_{t_2}\wedge \ldots\wedge e_{t_i}$ with $t_1<t_2<\cdots <t_i$.
The multidegree of $e_{t_1}\wedge e_{t_2}\wedge\cdots \wedge
e_{t_i}$ is the least common multiple of the elements
$u_{t_j}/[u_{t_j},u_p]$, $j=1,\ldots, i$.

The shifted complex
\[
\TT(-\multideg u_p)
\]
is a multigraded free resolution of $(S/L)(-\multideg u_p)$. We
denote the basis element of $T_i(-\multideg u_p)$ which
corresponds to $e_{t_1}\wedge e_{t_2}\wedge\cdots \wedge e_{t_i}$
by $b(p;\{t_1,\ldots,t_i\})$. Then $\multideg b(p;{t_1,\ldots,
t_i})= \multideg u_p + \multideg e_{t_1}\wedge e_{t_2}\wedge\cdots
\wedge e_{t_i}$, and  hence it is the least common multiple of
$u_p, u_{t_1},\ldots, u_{t_i}$.

The monomorphism $(S/L)(-\multideg u_p)\to S/H_\LL(q)$ induces a
comparison map
\[
\alpha\: \TT(-\multideg u_p)\To \FFF(q)
\]
of multigraded complexes. We let $\FFF(p)$ be the mapping cone of
$\alpha$. Then $\FFF(p)$ is a multigraded free $S$-resolution of
$H_\LL(p)$, and has the desired multigraded basis.

(b) (i) \implies (ii): Let $p\in \LL$ with $|N(p)|>1$, and let
$S\subset N(p)$ be a subset. By the definition of the differential
$\partial$ of $\FFF$ we have
\[
\partial b(p;S)=\sum_{q\in S} \pm v_q b(p;S\setminus\{q\}) +  \cdots
\]
where $v_q =\multideg b(p;S)/\multideg b(p;S\setminus\{q\})$.
Therefore the resolution can be minimal only if the multidegree of
$b(p;S\setminus\{q\})$ is a proper divisor the multidegree of
$b(p;S)$ for all $q$ in $S$.

By (a)
\[
\multideg  b(p;S)=x_{A}y_{B}\quad\text{and}\quad \multideg
b(p;S\setminus\{q\})=x_{A}y_{C},
\]
where $A=\ell(p)$, $B=\ell(p)^c\union\Union_{r\in S} \ell(r)^c$
and $C=\ell(p)^c\union\Union_{r\in S, r\neq q} \ell(r)^c$. Here,
for any subset $F\subset P$, we set $F^c=P\setminus F$.

It follows that $v_q=1$ if and only if $\Sect_{r\in S}\ell(r)=
\Sect_{r\in S, r\neq q}\ell(r)$. By Lemma \ref{better}(iii)  this
is equivalent to say that
\begin{eqnarray}
\label{not} \bigwedge \{r\: r\in S\}= \bigwedge  \{r\: r\in S,
r\neq q\}.
\end{eqnarray}
Hence if the resolution is minimal, then we do not have equality
in (\ref{not}) for any $S\subset N(p)$ and any $q\in S$. In
particular, for $S=N(p)$ we obtain the desired result.

(ii) \implies (i): Let $b(p;S)$ and $b(q;T)$ be two basis elements
with $|T|=|S|-1$. It suffices to show that in the following three
cases the coefficient of $b(q;T)$ in $\partial b(p;S)$ is either
$0$ or a monomial $\neq 1$:
\begin{itemize}
\item $p=q$ and $T\not\subset S$; \item $q<p$; \item $q\not < p$.
\end{itemize}
In the first case we show that $\multideg b(p;T)$ does not divide
$\multideg b(p;S)$. Otherwise we would have that $\Union_{r\in
T}\ell(r)^c\subseteq \Union_{r\in S}\ell(r)^c$. This would imply
that $\Sect_{r\in S}\ell(r)\subseteq \Sect_{r\in T}\ell(r)$, which
in turn would imply that $\bigwedge \{r\: r\in S\}\leq \bigwedge
\{r\: r\in T\}$. But then $\bigwedge \{r\: r\in N(p)\} =\bigwedge
\{r\: r\in N(p)\setminus(T\setminus S)\}$, a contradiction.

In the second case we have $\multideg b(p;S)=x_{\ell(p)}y_A$ and
$\multideg b(q;T)=x_{\ell(q)}y_B$ for some $A$ and $B$. If
$\multideg b(q;T)$ does not divide $\multideg b(p;S)$ then the
coefficient of $b(q;T)$ is $0$. Otherwise it is
$x_{\ell(p)\setminus \ell(q)}y_{A\setminus B}$. Since $q<p$ this
coefficient is not $1$.

In the last case $\ell(q)\not\subseteq \ell(p)$, and so $\multideg
b(q;T)$ does not divide $\multideg b(p;S)$. Hence the coefficient
of $b(q;T)$ is  $0$.
\end{proof}

\begin{Corollary}
\label{meet} If $\LL$ is a meet-distributive meet-semilattice,
then the finite multigraded free $S$-resolution given in Theorem
{\em \ref{free}} is minimal.
\end{Corollary}

\begin{proof}
By definition meet-distributive meet-semilattices  have the
property that for any element $p\in \LL$ the interval
$[\bigwedge\{q\: q\in N(p)\},p]$ is a Boolean lattice (of rank
$|N(p)|$). This implies condition (b)(ii) in Theorem \ref{free}.
\end{proof}

Note that condition (b)(ii) in Theorem \ref{free} is satisfied for
any meet-semilattice $\LL$ for which $|N(p)|\leq 2$ for all $p\in
\LL$. Other examples can easily be constructed, as follows: let
$\LL$ be a meet-semilattice satisfying the condition (b)(ii), and
let $p,q\in\LL$ such that $q\in N(p)$. Let $\LL'$ be the
meet-semilattice adding a new element $r$ with $q<r<p$. Then this
new meet-semilattice again satisfies (b)(ii).

An example of such a meet-semilattice is

\begin{center}
 \psset{unit=1.5cm}
\begin{pspicture}(0,-0.5)(2,2.5)
\psline(0.96,0.04)(0.54,0.46)
 \psline(0.5,0.54)(0.5,0.96)
 \psline(0.5,1.04)(0.5,1.46)
 \psline(0.54,1.54)(0.96,1.96)
 \psline(1,1.96)(1,1.54)
 \psline(0.96,1.46)(0.54,1.04)
 \psline(1,0.04)(1,0.46)
 \psline(1,0.54)(1,0.96)
 \psline(0.96,1.04)(0.54,1.46)
\psline(1.04,0.04)(1.46,0.46)
 \psline(1.5,0.54)(1.5,0.96)
 \psline(1.5,1.04)(1.5,1.46)
 \psline(1.46,1.04)(1.04,1.46)
 \psline(1.04,1.04)(1.46,1.46)
 \psline(1.46,1.54)(1.04,1.96)
 \rput(1,0){$\circ$}
 \rput(0.5,0.5){$\circ$}
 \rput(0.5,1){$\circ$}
 \rput(0.5,1.5){$\circ$}
 \rput(1,2){$\circ$}
 \rput(1,1.5){$\circ$}
 \rput(1,1){$\circ$}
 \rput(1,0.5){$\circ$}
 \rput(1.5,0.5){$\circ$}
 \rput(1.5,1){$\circ$}
 \rput(1.5,1.5){$\circ$}
  \rput(1,-0.4){${\mathcal L}$}
\end{pspicture}
\end{center}

Observe that $\LL$ is neither upper nor lower semimodular. The
resolution of $H_\LL$ is
\[
0\To S(-12)\To S^6(-10)\To S^9(-8)\dirsum S^6(-7)\To S^{11}(-6)\To
H_\LL\To 0.
\]

We close this section with discussing the regularity of $H_\LL$.
Recall that the regularity of a finitely generated graded
$S$-module $M$ is defined to be
\[
\reg M=\max\{j\: \beta_{i, i+j}(M)\neq 0 \text{ for some $i$}\}.
\]

\begin{Corollary}
\label{regularity} Let $\LL$ be a finite meet-semilattice and $P$
the poset of join irreducible elements in $\LL$. Then
\begin{enumerate}
\item[(a)] $\reg(H_\LL)\leq |P|+\max_{p\in\LL\atop  S\subset
N(p)}\bigl\{\deg p-\deg \bigwedge\{q\: q\in S\}-|S|\bigr\}$;
\item[(b)] if $\LL$ satisfies condition {\em (b)(ii)} in Theorem
{\em \ref{free}}, then
\[
\reg(H_\LL)= |P|+\max_{p\in\LL}\bigl\{\deg p-\deg \bigwedge\{q\:
q\in N(p)\}-|N(p)|\bigr\}.
\]
\end{enumerate}
\end{Corollary}

\begin{proof}
Since $\FFF$ is a possibly non-minimal free resolution of $H_\LL$
it follows that
\[
\reg H_\LL\leq \max\{\deg b(p;S)-|S|\}
\]
where the maximum is taken over all basis elements in the
resolution.

By our computation in the proof of Theorem \ref{free} one has
\[
\deg b(p;S)-|S|=|P|+\deg p-\deg \bigwedge\{q\: q\in S\}-|S|.
\]
This implies assertion (a).

If $\LL$ satisfies the condition (b)(ii) in Theorem \ref{free},
then our resolution is minimal and hence we have equality in
formula (a). Moreover, if $S'\subset S\subset N(p)$ with
$|S|=|S'|+1$,  then
\[
\deg \bigwedge\{q\: q\in S\}-\deg \bigwedge\{q\: q\in S'\}\geq 1.
\]
Hence
\[
\bigwedge\{q\: q\in S\}-\bigwedge\{q\: q\in N(p)\}\geq |N(p)|-|S|.
\]

\end{proof}

\section{The resolution of $H_\LL$ for a meet-distributive
meet-semilattice}

In this section we want to describe the differential $\partial$ in
the graded minimal free resolution $\FFF$ of $H_\LL$ when $\LL$ is
a meet-distributive meet-semilattice.

As we have seen in the previous section, a basis of $F_i$ is given
by the basis elements
\[
b(p;S),
\]
where $p\in \LL$ and $S\subset N(p)$ with $|S|=i$. Thus it amounts
to describe $\partial(b(p;S))$ for each such basis element. To
this end we introduce some notation:

Let $\LL$ be any meet-distributive meet-semilattice, and $P$ the
set of join-irreducible elements of $\LL$. We extend the partial
order on $P$ to a total order $<$.

For a subset $T\subset P$ and $q\in P$ we set
\[
\sigma(q;T)=|\{r\in T\: r<q\}|.
\]

For each  $q\in N(p)$, we have  $|\ell(p)\setminus\ell(q)|=1$. We
denote the unique element in $\ell(p)\setminus\ell(q)$ by
$p\setminus q$. Furthermore, for any subset $S\subset N(p)$ we set
$p\setminus S=\{p\setminus q\: q\in S\}$. With the notation
introduced we now have

\begin{Theorem}
\label{maps} For each $p\in \LL$ and each $S\subset N(p)$, one has
\[
\partial(b(p;S))=\sum_{q\in S}(-1)^{\sigma(p\setminus
q;p\setminus S)}(y_{p\setminus
q}b(p;S\setminus\{q\})-x_{p\setminus
q}b(q;q\wedge(S\setminus\{q\})).
\]
\end{Theorem}

Before we give the proof of the theorem we first note that
$q\wedge(S\setminus\{q\})\subset N(q)$ for all $q\in S$. This is
the case because by assumption $\LL$ is meet-distributive, so that
for any two distinct  lower neighbors $q_1$ and $q_2$ of $p$, the
element $q_1\wedge q_2$ is a lower neighbor of $q_1$ and $q_2$.

We also note that the differential defined in Theorem \ref{maps}
is multi-homogeneous. To see this, recall that $\multideg(b(p;S))$
is the least common multiple of $u_p$ and all $u_q$ with $q\in S$.
Since $u_q=y_{p\setminus q}u_p/x_{p\setminus q}$, we have $
\multideg(b(p;S\setminus \{q\}))=\multideg(b(p;S))/y_{p\setminus
q}, $ and $ \multideg(b(q;q\wedge (S\setminus
\{q\})))=\multideg(b(p;S))/x_{p\setminus q}. $ This shows that
$\partial$ is indeed multi-homogeneous.

\begin{proof}[Proof of \ref{maps}]
We use  the linear order $\prec$ on $\LL$ introduced  in the proof
of Theorem \ref{free}, and show by induction on  $p\in \LL$ that
the differential $\partial$ is given  on the free resolution
$\FFF(p)$ of $H_\LL(p)$ by the iterated mapping cone construction
as described in Theorem \ref{free}.

Recall that for $p\in \LL$ there is an exact sequence  of
multigraded $S$-modules
\[
0\To (S/L)(-\multideg u_p)\To S/H_\LL(q)\To S/H_\LL(p)\To 0,
\]
where $q\prec p$ is the element  in $\LL$ preceding $p$, and where
$L$ is the colon ideal
\[
H_\LL(q):u_p=(\{u_t/[u_t,u_p]\}_{t\in N(p)})=(y_{p\setminus t}\:
t\in N(p)).
\]
By induction hypothesis, the differential on $\FFF(q)$ is obtained
by iterated mapping cones from exact sequences as before.

Let $\CC=\TT(-\multideg u_p)$ be the shifted Taylor complex
associated with the sequence $y_{p\setminus t}$, $t\in N(P)$,
where the order of the sequence is given by the order of the
elements $p\setminus t$ in $P$. For a subset $S\in N(p)$,
$S=\{t_1,\ldots,t_i\}$ with $p\setminus t_1<p\setminus t_2<\cdots
<p\setminus t_i$, we denote the element $e_{t_1}\wedge
e_{t_2}\wedge\cdots\wedge  e_{t_i}\in T_i$ by $b(p;S)$.

Let $\alpha\: \CC\to \FFF(q)$ be a complex homomorphism extending
the map $$(S/L)(-\multideg u_p)\To S/H_\LL(q).$$ Then the
differential given by the mapping cone is defined as follows:
\[
\partial_i=(\partial^\TT_i+(-1)^i\alpha_i, \partial_{i+1}^{\FFF(q)}) \quad\text{for all}\quad i.
\]
Comparing this equation with the definition of $\partial$ in the
theorem it remains to  show that for each $S\subset N(p)$ we have:
\begin{enumerate}
\item[(i)] $\partial^{\TT}(b(p;S))=\sum_{q\in
S}(-1)^{\sigma(p\setminus q;p\setminus S)}y_{p\setminus
q}b(p;S\setminus\{q\})$, and

\item[(ii)]  $\alpha$ can be chosen such that
$$(-1)^i\alpha_i(b(p;S))=-\sum_{q\in S}(-1)^{\sigma(p\setminus
q;p\setminus S)}x_{p\setminus q}b(q;q\wedge(S\setminus\{q\}).$$
\end{enumerate}

Equation (i) is obvious, because this is exactly how the
differential in the Taylor complex is defined.

We conclude the proof of the theorem by showing that if  $\alpha$
is defined as in (ii), then $\alpha\: \CC\to \FFF(q)$ is a complex
homomorphism. This amounts to show that
\[
\partial^{\FFF(q)}_i\circ \alpha_i=\alpha_{i-1}\circ
\partial^{\TT}_i.
\]
To see this we choose $b(p;S)\in T_i$. Then
\begin{eqnarray}
\label{eq1}
 (\partial^{\FFF(q)}_i\circ
\alpha_i)(b(p;S))&=&(-1)^{i+1}\sum_{q\in S}(-1)^{\sigma(p\setminus
q;p\setminus S)}x_{p\setminus
q}\partial^{\FFF(q)}_i(b(q;q\wedge(S\setminus\{q\}))).
\end{eqnarray}
By our induction hypothesis we have that
\begin{eqnarray*}
\partial^{\FFF(q)}_i(b(q;q\wedge(S\setminus\{q\})))&=&
\sum_{q'\in S\setminus\{q\}}(-1)^{\sigma(p\setminus q';(p\setminus
S)\setminus\{p\setminus q\})}\bigl(y_{p\setminus q'}b(q;q\wedge
(S\setminus
\{q,q'\}))\\
&-&x_{p\setminus q'}b(q\wedge
q';q'\wedge[(q\wedge(S\setminus\{q\})\setminus\{q\wedge
q'\})])\bigr).
\end{eqnarray*}
Here we used that $q\setminus q\wedge q'=p\setminus q'$.

Substituting this in equation (\ref{eq1}) we get
\begin{eqnarray}
\label{eq2}
(\partial^{\FFF(q)}_i\circ \alpha_i)(b(p;S))=\hspace{9.2cm}\\
(-1)^{i+1} \sum_{q,q'\in S, q\neq q'}(-1)^{(\sigma(p\setminus
q;p\setminus S)+\sigma(p\setminus q';(p\setminus
S)\setminus\{p\setminus q\})}x_{p\setminus q}y_{p\setminus
q'}b(q;q\wedge (S\setminus \{q,q'\}).\nonumber
\end{eqnarray}
On the other hand
\begin{eqnarray}
\label{eq3}
 (\alpha_{i-1}\circ \partial^{\TT}_i)(b(p;S))=\sum_{q\in
S}(-1)^{\sigma(p\setminus q;p\setminus S)}y_{p\setminus
q}\alpha_{i-1}(b(p;S\setminus\{q\}))\hspace{2.5cm}\\
=(-1)^{i+1} \sum_{q,q'\in S, q\neq q'}(-1)^{(\sigma(p\setminus
q;p\setminus S)+\sigma(p\setminus q';(p\setminus
S)\setminus\{p\setminus q\})))}y_{p\setminus q}x_{p\setminus
q'}b(q';q'\wedge (S\setminus \{q,q'\}).\nonumber
\end{eqnarray}
Here we used that $q\setminus q\wedge q'=p\setminus q'$.

It follows that the right hand sides of the equations (\ref{eq2})
and (\ref{eq3}) coincide after exchanging $q$ and $q'$. This
concludes the proof.
\end{proof}

We would like to mention that our resolution is a cellular
resolution in the sense of Bayer and Sturmfels \cite{BS}, the
cells being cubes. Each basis element $b(p;S)$ can be identified
with the interval $[q,p]$ where $q$ is the meet of all elements in
$S$. Since $\LL$ is meet-distributive, this interval is a Boolean
lattice, and hence may be identified with a cube.

It would be desirable to have also an explicit description of the
differentials for the resolution of $H_\LL$ when $\LL$ is a
meet-irredundant meet-semilattice. Quite generally, according to
the iterated mapping cone construction described in Theorem
\ref{free}, the differentials in the resolution of $H_\LL$ for a
meet-irredundant meet-semilattice is of the form
\[
\partial(b(p;S))=\sum_{q\in S}(-1)^{\sigma(p\setminus
q;p\setminus S)}y_{p\setminus q}b(p;S\setminus\{q\})+\sum_{t\in
[r,p], t\neq p}c_tb(t;S_t),
\]
where
\begin{enumerate}
\item[(1)]  $r$ is the meet of all elements in $S$,

\item[(2)] $c_t=\lambda_t v_t$ with $\lambda_t\in K$ and $v_t$ the
monomial whose multidegree is
$\multideg(b(p;S))-\multideg(b(t;S_t))$,

\item[(3)]  $S_t$  is a set of lower neighbors of $t$ in the
interval $[r,p]$ with $|S_t|=|S|-1$.
\end{enumerate}

For example consider the following meet irredundant
meet-semilattice

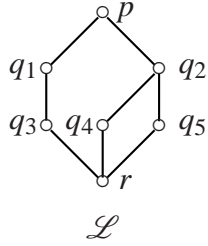
\begin{center}
 \psset{unit=1.5cm}
\begin{pspicture}(0,-0.5)(2,1.8)
\psline(0.96,0.04)(0.54,0.46)
 \psline(0.5,0.545)(0.5,0.96)
 \psline(0.54,1.04)(0.96,1.46)
 \psline(1.04,1.46)(1.46,1.04)
 \psline(1.5,0.96)(1.5,0.54)
 \psline(1.46,0.46)(1.04,0.04)
 \psline(1,0.04)(1,0.46)
 \psline(1.04,0.54)(1.46,0.96)
 \rput(1,0){$\circ$}
 \rput(0.5,0.5){$\circ$}
 \rput(0.5,1){$\circ$}
 \rput(1,1.5){$\circ$}
 \rput(1.5,1){$\circ$}
 \rput(1.5,0.5){$\circ$}
 \rput(1,0.5){$\circ$}
 \rput(1.2,0){$r$}
 \rput(0.3,0.5){$q_3$}
 \rput(0.3,1){$q_1$}
 \rput(1.2,1.5){$p$}
 \rput(1.8,1){$q_2$}
 \rput(1.8,0.5){$q_5$}
 \rput(0.8,0.5){$q_4$}
 \rput(1,-0.4){${\mathcal L}$}
\end{pspicture}
\end{center}
whose poset of join irreducible elements is
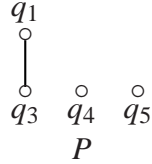
\begin{center}
 \psset{unit=1.5cm}
\begin{pspicture}(0,0)(2,1.5)
 \psline(0.5,0.545)(0.5,0.96)
 \rput(0.5,0.5){$\circ$}
 \rput(0.5,1){$\circ$}
 \rput(1.5,0.5){$\circ$}
 \rput(1,0.5){$\circ$}
 \rput(0.5,0.3){$q_3$}
 \rput(0.5,1.2){$q_1$}
 \rput(1.5,0.3){$q_5$}
 \rput(1,0.3){$q_4$}
 \rput(1,0){$P$}
\end{pspicture}
\end{center}

It is easy to see that in this case there are two, equally natural
choices, to define $\partial(b(p;\{q_1,q_2\})$, namely:
\begin{eqnarray*}
\partial(b(p;\{q_1,q_2\})=
-y_1y_3b(p;\{q_1\})+y_4y_5b(p;\{q_2\})-x_4x_5y_3b(q_1;\{q_3\})-x_1x_4x_5b(q_3;\{r\})\\
+x_1x_3y_4b(q_2;\{q_4\})+x_1x_3x_5b(q_4;\{r\}),\hspace{5.7cm}
\end{eqnarray*}
or,
\begin{eqnarray*}
\partial(b(p;\{q_1,q_2\})=
-y_1y_3b(p;\{q_1\})+y_4y_5b(p;\{q_2\})-x_4x_5y_3b(q_1;\{q_3\})-x_1x_4x_5b(q_3;\{r\})\\
+x_1x_3y_5b(q_2;\{q_5\})+x_1x_3x_4b(q_5;\{r\}).\hspace{5.7cm}
\end{eqnarray*}
Here we wrote for simplicity $x_i$ and $y_i$ instead of $x_{q_i}$
and $y_{q_i}$, respectively.

\section{On the Alexander dual of $H_\LL$}

For the convenience we introduce the following notation: let $I$
be a squarefree monomial ideal. Then  $I=I_\Delta$  for some
simplicial complex $\Delta$, and  we write $I^*$ for
$I_{\Delta^\vee}$. Here, as before, $\Delta^\vee$ is the Alexander
dual of the simplicial complex $\Delta$.

Let $\LL$ be a distributive lattice. In particular $\LL$ is a
poset and we may consider a poset ideal $\II\subset \LL$. Note
that any poset ideal $\II$ of $\LL$ is a (special)
meet-semilattice.

Let  $p\in \LL$, then the poset ideal
\[
\II_p=\{q\in \LL\: q\not\geq p\}
\]
is called {\em $1$-cogenerated}. It is clear that  for any poset
ideal $\II$ we have
\[
\II= \Sect_{p\in\LL\setminus\II}\II_p.
\]

We set $H_\II=(\{u_q\: q\in \II\}).$ Then

\begin{Lemma}
\label{easy} For any poset ideal $\II\in \LL$ we have
\[
H_\II=\Sect_{q\in\LL\setminus \II}H_{\II_q}\quad \text{and}\quad
H^*_\II=\sum_{q\in\LL\setminus \II}H^*_{\II_q}.
\]
\end{Lemma}

\begin{proof} In order to prove the first equation, it suffices to show that if
$\JJ$ and $\KK$ are two poset ideals in $\LL$, and $\II=\JJ\sect
\KK$, then $H_\II=H_\JJ\sect H_\KK$. It is clear that
$H_\II\subset H_\JJ\sect H_\KK$. Let $m\in H_\JJ\sect H_\KK$ a
monomial. Then there exist $p\in \JJ$ and $q\in \KK$ such that
$u_p|m$ and $u_q|m$. Let $t=p\wedge q$. Since $\LL$ is
distributive, we have $u_t=x_{\ell(p)\sect \ell(q)}y_{P\setminus
(\ell(p)\sect
\ell(q))}=x_{\ell(p)\sect\ell(q)}y_{(P\setminus\ell(p))\union
(P\setminus \ell(q))}$; hence $u_t|m$. Since $t\leq p$ and $t\leq
q$, it follows that $t\in \JJ\sect \KK=\II$. Therefore, $m\in
H_{\II}$.

Let $P$ be a monomial prime ideal. Then $\Sect_{q\in\LL\setminus
\II}H_{\II_q}\subset P$ if and only if $H_{\II_q}\subset P$ for
some $q$. Hence the assertion follows from (\ref{cover}).
\end{proof}

\begin{Theorem}
\label{alex} Let $\LL$ be a finite distributive lattice,
$P\subset\LL$ the poset of join irreducible elements of $\LL$, and
$\II\subset \LL$ a poset ideal of $\LL$. Then
\[
H^*_\II=(H^*_\LL,\{\prod_{r\in G(\ell(q))}y_r\: q\in\LL\setminus
\II\}),
\]
where $G(\ell(q))$ is the set of generators of the poset ideal
$\ell(q)\subset P$.
\end{Theorem}

\begin{proof}
By using Lemma \ref{easy} it suffices to prove the theorem for a
$1$-cogenerated poset ideal $\II_p$. In this case what we must
prove is
\[
H^*_{\II_p}=(H^*_\LL,\{\prod_{r\in G(\ell(q))}y_r\: q\geq p\}).
\]

Let $x_Ay_B$ be a squarefree monomial with $A,B\subset P$. Then
$x_Ay_B \in H^*_{\II_p}$ if and only if  $A\sect \ell(r)\neq
\emptyset$, or $B\sect \ell(r)^c\neq \emptyset$ for all $r\not\geq
p$.

Let $T=(H^*_\LL,\{\prod_{r\in G(\ell(q))}y_r\: q\geq p\})$. We
first show that $T\subset H^*_{\II_p}$. Since $H_{\II_p}\subset
H_\LL$ it follows that $H_\LL^*\subset H^*_{\II_p}$. Moreover,
suppose that for some  $q\geq p$ the monomial $\prod_{r\in
G(\ell(q))}y_r$ does not belong to $H^*_{\II_p}$. Then there
exists $t\not\geq p$ such that $G(\ell(q))\sect
\ell(t)^c=\emptyset$, equivalently $G(\ell(q))\subset \ell(t)$.
Hence $\ell(q)\subset \ell(t)$. However, since $q\geq p$, we have
$\ell(p)\subset \ell(q)$, so that $\ell(p)\subset \ell(q)$, a
contradiction.

It remains to show that $H^*_{\II_p}\subset T$.

Suppose $B=\emptyset$. Then $A\sect \ell(\hat{0})=\emptyset$ since
$\ell(\hat{0})=\emptyset$ and also $B\sect
\ell(\hat{0})^c=\emptyset$, a contradiction.

Suppose $A = \emptyset$.  Let $\Delta^\vee$ denote the simplicial
complex whose Stanley--Reisner ideal is equal to $H^*_{\II_p}$ and
$\Delta^\vee_y$ the restriction of $\Delta^\vee$ on the vertex set
$\{y_t \: t \in P \}$.  Then the facets of $\Delta^\vee_y$ are
$\{y_t \: t \in \II \}$, where $\II$ is a maximal poset ideal of
$P$ which does not contain $\ell(p)$.  Such a poset ideal is of
the form $P \setminus \{ t \in P \: t \geq h \}$ with $h \in
G(\ell(p))$.  If $y_B$ belongs to $H^*_{\II_p}$, then $B$ is
contained in no facet of $\Delta^\vee_y$.  Hence, for each $h \in
G(\ell(p))$, there is $h' \in P$ with $h' \geq h$ such that $h'
\in B$.  Let $\II_0$ denote the poset ideal of $P$ consisting of
all $t \in P$ with $t \leq h'$ for some $h \in G(\ell(p))$.  Let
$q \in \LL$ with $\ell(q) = \II_0$.  It then follows that
$\prod_{r\in G(\ell(q))}y_r$ divides $y_B$.

Finally we consider the case that $A\neq\emptyset$, and
$y_B\not\in H_{\II_p}^*$. We will show that in this case
$x_Ay_B\in H_\LL^*$. In fact, since $y_B\not\in H_{\II_p}^*$,
there exists $r\not\geq p$ such that $B\sect \ell(r)^c=\emptyset$,
equivalently $B\subset \ell(r)$. Let $(B)\subset P$ be the poset
ideal generated by $B$. Then there exists $t\in\LL$ such that
$\ell(t)=(B)$. Since $\ell(t)=(B)\subset\ell(r)$ it follows that
$t\leq r$, and hence $t\in \II_p$.

Suppose $x_Ay_B\not\in H_\LL^*$, then $a\not\leq b$ for all $a\in
A$ and $b\in B$. This implies that $A\sect (B)=A\sect
\ell(t)=\emptyset$. This is  a contradiction because also $B\sect
\ell(t)^c=\emptyset$.
\end{proof}

Recall from \cite[Theorem 2.4]{HH} that if $G$ is a
Cohen--Macaulay bipartite graph on the vertex set $V\union V'$
with $V\sect V'=\emptyset$ and $|V|=|V'|$, then there exists a
partial order $<$ on $V$ such that the distributive lattice
$\JJ(P)$ with $P = (V, <)$ satisfies $H_{\JJ(P)}^*=I(G)$.  We
write $\LL(G)$ for the distributive lattice $\JJ(P)$.

\begin{Theorem}
\label{takayuki} Let $\Delta$ be a simplicial complex on the
vertex set $V\union V'$ with $V\sect V'=\emptyset$  and
$|V|=|V'|$. Suppose that
\begin{enumerate}
\item[(1)] there is no $F\in \mathcal{F}(\Delta)$ with $F\subset
V$,

\item[(2)] $G=\{F\in\mathcal{F}(\Delta)\: F\sect V\neq
\emptyset,\quad F\sect V'\neq \emptyset\}$ is a Cohen--Macaulay
bipartite graph with no isolated vertex.
\end{enumerate}
Then the following conditions are equivalent:

\begin{enumerate}
\item[(a)] $S/I(\Delta)$ is Cohen--Macaulay;

\item[(b)] The simplicial complex $\Gamma$ with
$I_\Gamma=I(\Delta)$ is pure;

\item[(c)] There exists a poset ideal $\II\subset \LL(G)$
containing all join-irreducible elements of $\LL(G)$ such that
$H_\II^*=I(\Delta)$.
\end{enumerate}
\end{Theorem}

The following pictures show  examples of simplicial complexes
satisfying the conditions (1) and (2) of Theorem \ref{takayuki}.

\begin{center}
 \psset{unit=1.5cm}
\begin{pspicture}(0,0)(5,2)
 \psline(0.5,0.44)(0.5,1.36)
 \psline(1,0.79)(1,1.71)
 \psline(1.5,0.79)(1.5,1.71)
 \psline(2,0.79)(2,1.71)
 \psline(0.54,1.36)(1.46,0.79)
 \psline(1.04,1.71)(1.96,0.79)
 \psline(1.54,0.75)(1.96,0.75)
 \psline(0.54,0.44)(1.96,0.71)
 \rput(0.5,0.4){$\circ$}
 \rput(0.5,1.4){$\circ$}
 \rput(1,0.75){$\circ$}
 \rput(1,1.75){$\circ$}
 \rput(1.5,0.75){$\circ$}
 \rput(1.5,1.75){$\circ$}
 \rput(2,0.75){$\circ$}
 \rput(2,1.75){$\circ$}
 \rput(1.25,0){$\Delta$}

  \psline(3.5,0.44)(3.5,1.36)
 \psline(4,0.79)(4,1.71)
 \psline(4.5,0.79)(4.5,1.71)
 \psline(5,0.79)(5,1.71)
 \psline(3.54,1.36)(4.46,0.79)
 \psline(4.04,1.71)(4.96,0.79)
 \pspolygon[style=fyp, fillcolor=medium](3.5,0.4)(5,0.75)(4.5,0.75)
 \rput(3.5,0.4){$\circ$}
 \rput(3.5,1.4){$\circ$}
 \rput(4,0.75){$\circ$}
 \rput(4,1.75){$\circ$}
 \rput(4.5,0.75){$\circ$}
 \rput(4.5,1.75){$\circ$}
 \rput(5,0.75){$\circ$}
 \rput(5,1.75){$\circ$}
 \rput(4.25,0){$\Delta'$}
\end{pspicture}
\end{center}

The facet ideal of the simplicial complex $\Delta$ is
Cohen--Macaulay, and that of $\Delta'$ is not Cohen-Macaulay. In
fact, the distributive lattice $\LL$ and its poset $P$ of join
irreducible elements corresponding to the bipartite graph in
$\Delta$ and $\Delta'$ is in both cases
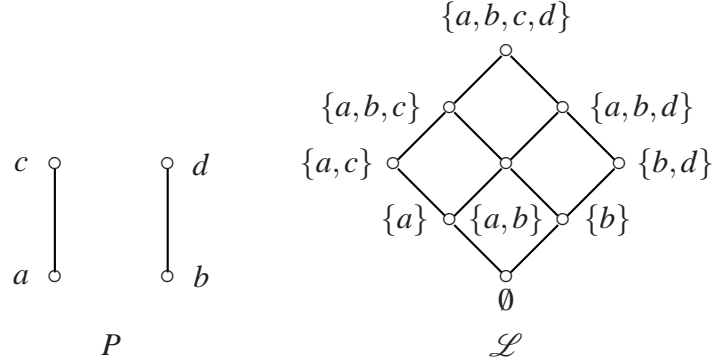
\begin{center}
 \psset{unit=1.5cm}
\begin{pspicture}(0,-1)(7,3)
\psline(1,0.04)(1,0.96)
 \psline(2,0.04)(2,0.96)
 \rput(1,0){$\circ$}
 \rput(1,1){$\circ$}
 \rput(2,0){$\circ$}
 \rput(2,1){$\circ$}
 \rput(0.7,0){$a$}
 \rput(0.7,1){$c$}
 \rput(2.3,0){$b$}
 \rput(2.3,1){$d$}
 \rput(1.5,-0.6){$P$}

\psline(4.96,0.04)(4.54,0.46)
 \psline(4.54,0.54)(4.96,0.96)
 \psline(5.04,0.96)(5.46,0.54)
 \psline(5.46,0.46)(5.04,0.04)
 \psline(5.54,0.54)(5.96,0.96)
 \psline(4.46,0.54)(4.04,0.96)
 \psline(4.04,1.04)(4.46,1.46)
 \psline(4.54,1.46)(4.96,1.04)
 \psline(5.04,1.04)(5.46,1.46)
 \psline(5.54,1.46)(5.96,1.04)
 \psline(4.54,1.54)(4.96,1.96)
 \psline(5.04,1.96)(5.46,1.54)
 \rput(5,0){$\circ$}
 \rput(4.5,0.5){$\circ$}
 \rput(4,1){$\circ$}
 \rput(5,1){$\circ$}
 \rput(5.5,0.5){$\circ$}
 \rput(6,1){$\circ$}
 \rput(4.5,1.5){$\circ$}
 \rput(5.5,1.5){$\circ$}
 \rput(5,2){$\circ$}
 \rput(5,-0.2){$\emptyset$}
 \rput(4.1,0.5){$\{a\}$}
 \rput(3.5,1){$\{a,c\}$}
 \rput(5,0.5){$\{a,b\}$}
 \rput(5.9,0.5){$\{b\}$}
 \rput(6.5,1){$\{b,d\}$}
 \rput(3.8,1.5){$\{a,b,c\}$}
 \rput(6.2,1.5){$\{a,b,d\}$}
 \rput(5,2.3){$\{a,b,c,d\}$}
 \rput(5,-0.6){$\LL$}
\end{pspicture}
\end{center}

The simplicial complex $\Delta$ corresponds to the ideal
$$\II=\{\emptyset, \{a\}, \{b\},
\{a,b\},\{a,c\},\{b,d\},\{a,b,c\}\}.$$

Since all poset ideals of $\LL$ are generated by at most two
elements, it follows from Theorem \ref{alex} that the simplicial
complex $\Delta'$ cannot correspond to any poset ideal in $\LL$.
Therefore, by Theorem \ref{takayuki} it cannot be Cohen-Macaulay.

\begin{proof}[Proof of Theorem \ref{takayuki}.]
Since every Cohen--Macaulay simplicial complex is pure, one has
(a) $\Rightarrow$ (b).  Moreover, since Theorem \ref{main}
guarantees that $H_\II$ has a linear resolution, it follows from
Lemma \ref{duality} that (c) $\Rightarrow$ (a).

We now prove that (b) $\Rightarrow$ (c).  Let $V = \{ x_1, \ldots,
x_n \}$ and $V' = \{ y_1, \ldots, y_n \}$.  Since $\Gamma$ is pure
and since $V$ is a facet of $\Gamma$, it follows that each facet
of $\Gamma$ is a facet of $\Gamma_0$, where $\Gamma_0$ is a
simplicial complex on $V \union V'$ with $I_{\Gamma_0} = I(G)$. In
other words, each minimal nonface of $\Gamma^\vee$ is a minimal
nonface of $\Gamma_0^\vee$. Thus we may regards that the minimal
set $\II^\flat$ of monomial generators of $I_{\Gamma^\vee}$ is a
subset of $\LL(G)$. Now, what we must prove is that $\II^\flat$ is
a poset ideal of $\LL(G) = \JJ(P)$, where $P = (V, <)$ is the
poset consisting of all join-irreducible elements of $\LL(G)$.
Suppose, on the contrary, that $\II^\flat$ is not a poset ideal,
and choose two elements $\delta$ and $\xi$ of $\LL(G)$ with
$\delta \in \II^\flat$ and $\xi \not\in \II^\flat$ such that
$\delta$ covers $\xi$ in $\LL(G)$.  To simplify the notation, we
will assume that $\delta = \{ x_1, \ldots, x_k \}$ and $\xi = \{
x_1, \ldots, x_{k-1} \}$.  Thus $\{ y_1, \ldots, y_k, x_{k+1},
\ldots, x_n \}$ is a facet of $\Gamma$ and $\{ y_1, \ldots,
y_{k-1}, x_k, x_{k+1}, \ldots, x_n \}$ is not a facet of $\Gamma$.
Thus there is a monomial generator $u$ of $I(\Delta)$ which
divides $y_1 \cdots y_{k-1} x_k x_{k+1} \cdots x_n$.  However,
since $\{ y_1, \ldots, y_{k-1}, x_{k+1}, \ldots, x_n \}$ is a face
of $\Gamma$, it follows that the variable $x_k$ must appear in the
support of $u$.  Hence $u = x_k y_j$ with $1 \leq j \leq k - 1$.
Then \cite[Theorem 3.4]{HH} says that $x_k < x_j$ in $P$.  This is
impossible, since $\xi$ is a poset ideal of $\LL(G)$.
Consequently, it turns out that $\II^\flat$ is a poset ideal of
$\LL(G)$.

Finally, in case that $\II^\flat$ does not contain of a
join-irreducible element $x_i$ of $\LL(G)$, the vertex $y_i$
belongs to all facets of $\Gamma$.  This is impossible, since $G$
possesses no isolated vertex.  This completes the proof of (b)
$\Rightarrow$ (c).
\end{proof}

\begin{Corollary}
\label{grafted} Let $\Delta$ be a simplicial complex on the vertex
set $V=\{v_1,\ldots, v_n\}$, and let $W=\{w_1,\ldots, w_n\}$ be a
vertex set with $W\sect V=\emptyset$. Let $\Gamma$ be the
simplicial complex on the vertex set $V\union W$ whose facets are
those of $\Delta$ and all the edges $\{v_i,w_i\}$ for $i=1,\ldots,
n$. Then the facet ideal of $\Gamma$ is Cohen-Macaulay.
\end{Corollary}

\begin{proof}
Our work is to show that the simplicial complex $\Sigma$ with
$I_{\Sigma} = I(\Gamma)$ is pure.  Let $F = \{ v_i \: i \in A \}
\union \{ w_j \: j \in B \}$ be a face of $\Sigma$;  then  $A \cap
B = \emptyset$.  If $A \union B \neq [n]$, then $F \cup \{ w_i \:
i \in [n] \setminus (A \union B) \}$ is a face of $\Sigma$. Thus
all facets of $\Sigma$ have the cardinality $n$.  Hence $\Sigma$
is pure, as desired.
\end{proof}

The results of Theorem \ref{main} and Theorem \ref{alex} can be
extended as follows. Let $P$ be a poset. Recall that a {\em poset
coideal} of $P$ is a subset $J\subset P$ with the property that
for each $p\in J$ and each $q\in P$ with $q\geq p$ one has $q\in
J$. The minimal elements in $J$ are called the {\em cogenerators}.
The set of cogenerators of $J$ will be denoted by $G(J)$.

Now  let $\LL$ be a finite distributive lattice, and let
$\II\subset \LL$ be a poset ideal, and $\JJ$ a poset coideal  in
$\LL$. Then $H_\II$ and $H_\JJ$ have linear resolutions. We know
this for $H_\II$ by Theorem \ref{main} and for $H_\JJ$ it follows
by the same theorem using the fact that the dual of $\LL$ (where
the order of the elements of $\LL$ is just reversed) is again a
distributive lattice. What can be said about $H_\II\sect H_\JJ$?
The reader might expect that this ideal has  again a linear
resolution. However this is not the case. For example, consider
the Boolean lattice ${\mathcal B}_3$ of rank 3, and let
$\II={\mathcal B}_3\setminus\{\hat{1}\}$ and $\JJ={\mathcal
B}_3\setminus\{\hat{0}\}$. Then $H_\II\sect H_\JJ$ does not have a
linear resolution.

However in the positive direction we have

\begin{Proposition} \label{positive}
Let $\II$ be a poset ideal and $\JJ$ a poset coideal in $\LL$.
Then

\begin{enumerate}

\item[(a)] $\rank \LL\leq  \reg(H_\II\sect H_\JJ)\leq \rank\LL
+1$, if $\LL=\II\union\JJ$.

\item[(b)] $(H_\II\sect H_\JJ)^*=(H_\LL^*, \{\prod_{r\in
G(\ell(q))}y_r\: q\in\LL\setminus \II\},\{\prod_{r\in
G(\ell(q)^c)}x_r\: q\in\LL\setminus \JJ\})$.
\end{enumerate}
\end{Proposition}

\begin{proof}
(a) Consider the long exact Tor-sequence
\[
\cdots\to \Tor_{i+1}(K,H_\II+H_\JJ)\to \Tor_i(K,H_\II\sect
H_\JJ)\to \Tor_i(K,H_\II)\dirsum \Tor_i(K,H_\JJ)\to \cdots
\]
arising from the short exact sequence
\[
0\To H_\II\sect H_\JJ\To H_\II\dirsum H_\JJ\To  H_\II+H_\JJ\To 0.
\]
Since $H_\LL=H_\II+H_\JJ$, the ideals $H_\II$, $H_\JJ$ and
$H_\II+H_\JJ$ have a linear resolution by Theorem \ref{main}. It
follows that $\Tor_i(K,H_\II)_j=\Tor_i(K,H_\JJ)_j=0$ for $j\neq
i+\rank \LL$, and $\Tor_{i+1}(K,H_\II+H_\JJ)_j=0$ for $j\neq
i+1+\rank \LL$. Thus the assertion follows from the long exact
Tor-sequence.

(b) Since $(H_\II\sect H_\JJ)^*=H_\II^*+H_\JJ^*$, the claim
follows Theorem \ref{alex}.
\end{proof}

Consider the following example.
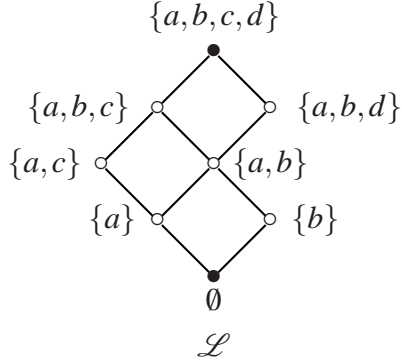
\begin{center}
 \psset{unit=1.5cm}
\begin{pspicture}(3,-1)(7,3)
\psline(4.96,0.04)(4.54,0.46)
 \psline(4.54,0.54)(4.96,0.96)
 \psline(5.04,0.96)(5.46,0.54)
 \psline(5.46,0.46)(5.04,0.04)
 \psline(4.46,0.54)(4.04,0.96)
 \psline(4.04,1.04)(4.46,1.46)
 \psline(4.54,1.46)(4.96,1.04)
 \psline(5.04,1.04)(5.46,1.46)
 \psline(4.54,1.54)(4.96,1.96)
 \psline(5.04,1.96)(5.46,1.54)
 \rput(5,0){$\bullet$}
 \rput(4.5,0.5){$\circ$}
 \rput(4,1){$\circ$}
 \rput(5,1){$\circ$}
 \rput(5.5,0.5){$\circ$}
 \rput(4.5,1.5){$\circ$}
 \rput(5.5,1.5){$\circ$}
 \rput(5,2){$\bullet$}
 \rput(5,-0.2){$\emptyset$}
 \rput(4.1,0.5){$\{a\}$}
 \rput(3.5,1){$\{a,c\}$}
 \rput(5.5,1){$\{a,b\}$}
 \rput(5.9,0.5){$\{b\}$}
 \rput(3.8,1.5){$\{a,b,c\}$}
 \rput(6.2,1.5){$\{a,b,d\}$}
 \rput(5,2.3){$\{a,b,c,d\}$}
 \rput(5,-0.6){$\LL$}
\end{pspicture}
\end{center}
Here we take in $\LL$ the poset ideal $$\II=\{\emptyset,
\{a\},\{b\}, \{a,b\}, \{a,c\}, \{a,b,c\}, \{a,b,d\}\},$$ and the
poset coideal
$$\JJ= \{ \{a\},\{b\}, \{a,b\}, \{a,c\}, \{a,b,c\},
\{a,b,d\},\{a,b,c,d\}\}.$$

Then $H_\II\sect H_\JJ=(avwx, buwx, acvx, abwx, abcx, abdw)$. Thus
this intersection is generated by all generators of $H_\LL$ except
$u_{\hat{0}}$ and $u_{\hat{1}}$, as indicated in the picture. The
resolution of $H_\II\sect H_\JJ$ is linear, namely
\[
0\To S(-6)\To S(-5)^6\To S(-4)^6\To H_\II\sect H_\JJ\To 0.
\]

Quite generally it would be interesting to know when  $H_\II\sect
H_\JJ= H_{\II\sect \JJ}$, and when an ideal of the  form
$H_{\II\sect \JJ}$ has a linear resolution. Of particular interest
are the following cases:
\begin{enumerate}
\item[(1)] $H=(\{u_p\}_{p\in \LL\setminus \{\hat{0},\hat{1}\}})$;
\item[(2)] $H=(\{u_p\: r\leq \rank p\leq s\})$ for some $r$ and
$s$ with $0<r\leq s<\rank \LL$.
\end{enumerate}

\end{document}